\documentclass[hidelinks,onefignum,onetabnum]{siamart220329}

\usepackage{lipsum}
\usepackage{amsfonts}
\usepackage{graphicx}
\usepackage{epstopdf}

\usepackage{algorithm, algpseudocode}
\algrenewcommand\algorithmicrequire{\textbf{Input:}}
\algrenewcommand\algorithmicensure{\textbf{Output:}}

\usepackage{amsmath}
\usepackage{amssymb}
\usepackage{bm}
\usepackage{bbm}

\ifpdf
  \DeclareGraphicsExtensions{.eps,.pdf,.png,.jpg}
\else
  \DeclareGraphicsExtensions{.eps}
\fi

\newcommand{\lf}{h} 
\newcommand{\sympmap}{\bm{F}} 
\newcommand{\Ginv}{G_{\mathrm{Inv}}}
\newcommand{\Einv}{E_{\mathrm{Inv}}}

\newcommand{\wbd}{w_{\mathrm{bd}}}
\newcommand{\Wbd}{W_{\mathrm{bd}}}
\newcommand{\Ebd}{E_{\mathrm{bd}}}
\newcommand{\hbd}{h_{\mathrm{bd}}}
\newcommand{\KZZ}{K_{ZZ}}
\newcommand{\EK}{E_K}
\newcommand{\WB}[1]{\mathcal{WB}[#1]}

\newcommand{\Rbb}{\mathbb{R}}

\newcommand{\Tbb}{\mathbb{T}}

\newcommand{\abs}[1]{\left | #1 \right |}
\newcommand{\norm}[1]{\left \lVert #1 \right \rVert}
\newcommand{\ip}[2]{\left\langle #1 , \ #2 \right\rangle }

\newcommand{\dif}{\hspace{1pt} \mathrm{d}}

\newsiamremark{remark}{Remark}

\headers{Level Set Learning for Poincar\'e Plots of Symplectic Maps}{M. Ruth and D. Bindel}

\title{Level Set Learning for Poincar\'e Plots of Symplectic Maps
    \thanks{
        Submitted to the editors November 30, 2023.
        \funding{This work was supported by the National Science Foundation Graduate Research Fellowship under Grant No. DGE-1650441 and by the and by a grant from the Simons Foundation (No. 560651, D B).}
    }
}

\author{Maximilian Ruth\thanks{Center for Applied Mathematics, Cornell University, Ithaca, NY
  (\email{mer335@cornell.edu})}
\and David Bindel\thanks{Department of Computer Science, Cornell University, Ithaca, NY }
}

\usepackage{amsopn}

\ifpdf
\hypersetup{
  pdftitle={Ruth-Bindel-KernelLabel},
  pdfauthor={M. Ruth, D. Bindel}
}
\fi

\begin{document}

\maketitle

\begin{abstract}
Many important qualities of plasma confinement devices can be determined via the Poincar\'e plot of a symplectic return map. 
These qualities include the locations of periodic orbits, magnetic islands, and chaotic regions of phase space. 
However, every evaluation of the magnetic return map requires solving an ODE, meaning a detailed Poincar\'e plot can be expensive to create. 
Here, we propose a kernel-based method of learning a single labeling function that is approximately invariant under the symplectic map. 
From the labeling function, we can recover the locations of invariant circles, islands, and chaos with few evaluations of the underlying symplectic map.
Additionally, the labeling function comes with a residual, which serves as a measure of how invariant the label function is, and therefore as an indirect measure of chaos and map complexity.
\end{abstract}



\section{Introduction}

\begin{figure}
    \centering
    \includegraphics[width=\textwidth]{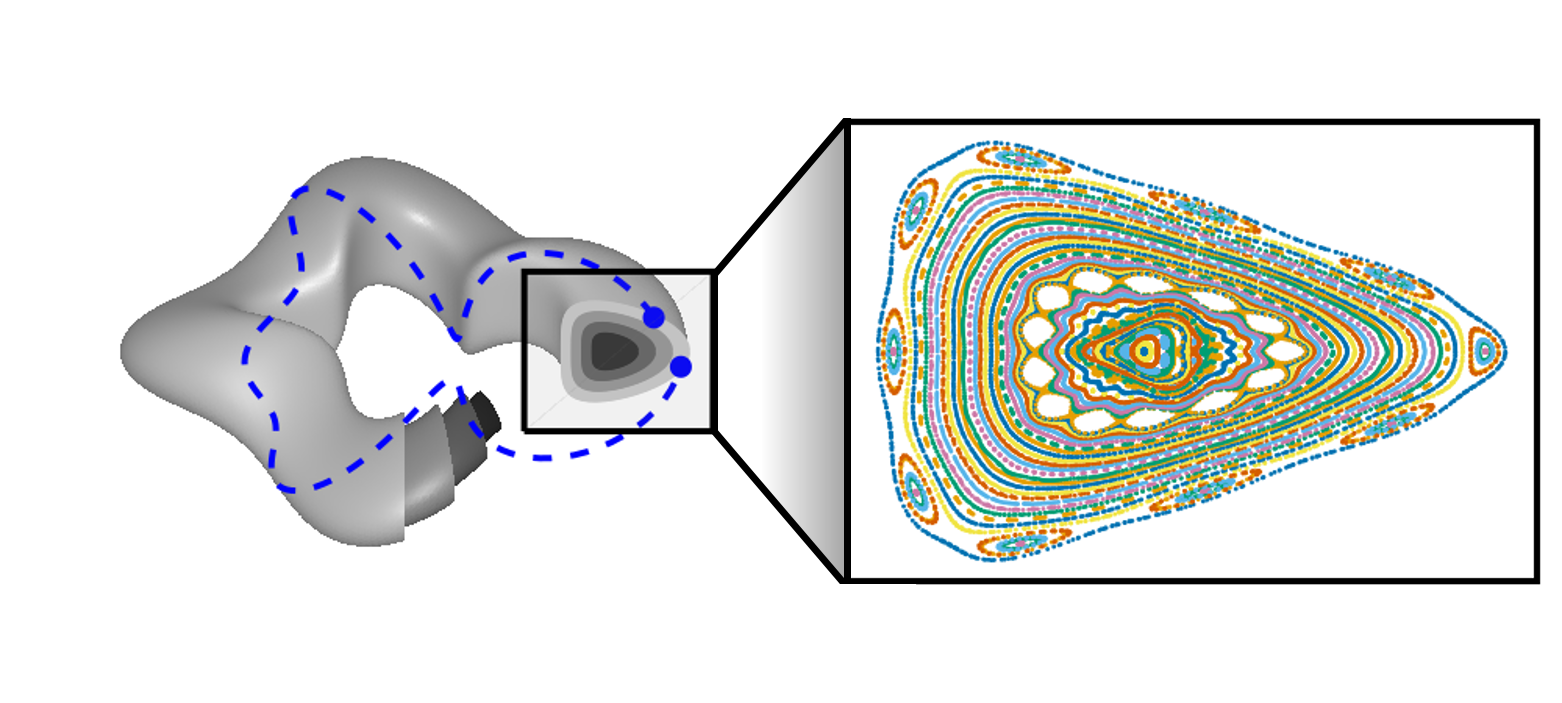}
    \caption{(left) A schematic visualization of a stellarator. The nested surfaces are tangent to the magnetic field, so particles are trapped to leading order. An example magnetic field line is given in blue, where the symplectic return map is defined by following the field line back to a given cross-section. (right) An example Poincar\'e plot, explored further in \cref{sec:stellarator-example}. The plot consists primarily of nested invariant circles. The outer ring has an island chain, and the center has a small amount of chaos.}
    \label{fig:cross-section}
\end{figure}

Magnetic field line flow is a useful first approximation to particle trajectories in stellarators \cite{hazeltine2003, helander2014, imbert-gerard2020}, a type of toroidal magnetic confinement device for plasmas.
When the magnetic flow is tangent to a set of nested surfaces, we expect good confinement of particles.
This can be visualized via a cross-section of the magnetic flow (see \cref{fig:cross-section}), where nested surfaces correspond to nested circles of a magnetic symplectic return map. 
Other invariant structures on the cross-section, such as islands and chaos, usually degrade confinement.
As such, when designing for good stellarator configurations, good methods to \textit{quickly} and \textit{automatically} analyze the invariant structures of a map are required.

There are many existing methods for finding invariant structures.
The simplest method -- the Poincar\'e plot -- is useful, but requires human intervention to determine the locations of islands and chaos.
More complicated visualization methods via ergodic averaging \cite{mezic1999} and clustering \cite{levnajic2010} improve upon the straightforward Poincar\'e plot, but still require human intervention.
Building upon the ergodic averaging work, it has been shown that the convergence rate of a weighted Birkhoff average \cite{Das2017, sander2020} can be used to quickly differentiate between chaos and invariant circles. However, this still does not easily differentiate islands from the core.

Another group of methods finds invariant and approximately invariant structures directly.
The most direct method is the parameterization method \cite{haro2016}, which uses a Newton iteration on the parameterization used in the Kolmogorov–Arnold–Moser (KAM) theorem.
Invariant circles found via the parameterization can be proven to be invariant via a numerical KAM theorem \cite{Figueras2017}.
However, good initial guesses for the parameterization method are not always available.
A comparable method in stellarator physics is quadratic flux minimization (QFMin) in maps \cite{dewar1992} and 3D divergence-free fields \cite{dewar1994}. 
QFMin is a more robust alternative to the parameterization method, but it comes at the cost of averaging over structures such as islands and chaos.
Recently, it has been proposed to solve an elliptic temperature diffusion PDE with boundary conditions to determine whether good nested invariant circles exist \cite{paul2022}, which has been proven to recover approximately invariant circles even under structure-breaking perturbations \cite{drivas2022}. 

A final set of methods is to first approximate the Perron-Frobenius operator or its dual the Koopman operator, and then to find the spectrum of that operator. 
The Perron-Frobenius based method uses an Ulam discretization to find approximately invariant distributions of particles \cite{froyland2009, froyland2014}. 
This can be viewed as a member of the broader problem of Lagrangian coherent structure detection \cite{hadjighasem2017}. 
There is a large amount of recent research on the Koopman operator \cite{kutz2016dynamic}, with recent work specifically focused on symplectic maps in \cite{govindarajan2019, colbrook2022}.

Many of the above methods have a similar problem: they require many iterations of the symplectic map to converge. 
Every symplectic map evaluation for a stellarator requires solving an ODE on the field line, so the speeds of these methods are often limited by the cost of evaluating the map. 
Additionally, the stellarator optimization process requires analyzing many configurations (and hence many maps), amplifying this inefficiency.
One way to ameliorate the evaluation cost is to approximate the symplectic map.
The simplest option is linear approximation (e.g.~using Chebyshev polynomials \cite{boyd2001}), but linear approximations break the symplecticity of the map, which can harm the performance of the above methods --- especially those using long trajectories, which are sensitive to structural perturbations.
Alternatively, many specialized symplectic approximation methods have been proposed, including the H\'enonNet \cite{Burby2021}, SympNet \cite{Jin2020}, and a symplectic Gaussian process approach \cite{Rath2021}. 
However, these methods have the shared drawback that good symplectic approximation requires nonlinear approximation, meaning a complicated and potentially lengthy training phase must be undertaken.

In this paper, we propose an alternative sample-efficient way of determining the invariant sets of a symplectic map.
Specifically, we find a smooth nonconstant ``label'' function that is approximately invariant under composition with the symplectic map (i.e.~it is approximately an eigenvector of the Koopman operator with eigenvalue $1$).
When the map actually consists of nested invariant circles, their locations are approximately recovered as the level sets of the invariant function. 
The fit of the smooth function also comes with a residual, which serves as a measure of how invariant the label function is. 
This residual is smooth with respect to the symplectic map and does not require \textit{a priori} knowledge of island or chaos location, making it appealing as a potential optimization target.

We present two methods for finding the invariant function.
The first is a boundary value problem, where we assume the user approximately knows the locations of inner and outer invariant circles, and would like to determine whether there are invariant circles foliating the region. 
The boundary values force the invariant function to be nonzero and different on the boundaries, so the existence of a good invariant function requires that no trajectory quickly reaches on boundary from the other.
This is posed as a linear least-squares problem, where the function is assumed to be in a reproducing kernel Hilbert space.
The second method is an eigenvalue problem, where we assume the user only knows the outer edge of a domain of interest.
In this method, boundary conditions force the function to be zero on the edge, and the eigenvalue problem forces the invariant function to be smoothly nonzero on the interior. 
The problem is posed as a Rayleigh quotient, and the eigenvalue gives the measure of fit.
Under typical conditions, the cost of both methods is dominated by time to evaluate the symplectic map, and the number of map evaluations is less than around 1000. We show that this number of evaluations is plenty for qualitative analysis, making both methods fast enough to use for stellarator optimization. Code for these methods is implemented in the \texttt{SymplecticMapTools.jl}\footnote{\href{https://github.com/maxeruth/SymplecticMapTools.jl}{https://github.com/maxeruth/SymplecticMapTools.jl}} Julia package.

In \cref{sec:description-of-method} we describe our approach to finding smooth, invariant, and nonzero functions of symplectic maps.
The section begins with an introduction to reproducing kernel Hilbert spaces and our measure for enforcing invariance.
Then, we present our two methods for finding nonzero functions. 
In the following sections, we analyze two detailed examples of our methods. 
In the first example (\cref{sec:stellarator-example}), we solve the eigenvalue problem for on an optimized stellarator configuration.
We measure the invariance of the function via comparison against a weighted Birkhoff average, and demonstrate how the invariance changes with the number of samples and kernel width.
We find the label function can resolve more detailed invariant structures as the number of samples increases and the kernel width decreases.
In the second example (\cref{sec:standard-map}), we solve the boundary value problem on the standard map. 
In particular, we analyze the boundary value problem in the transition to chaos, showing how the least squares residual measures this process.
Finally, we conclude in \cref{sec:conclusions}.

\section{Description of the Method}
\label{sec:description-of-method}

As outlined in the introduction, the goal of this paper is to find \textit{approximately invariant}, \textit{smooth}, and \textit{nonconstant} functions.
In this section, we describe how each of these qualities is achieved. 
We begin in \cref{subsec:preliminaries} with some preliminaries, where we define symplecticity and invariance, and introduce the Birkhoff average that will be used in \cref{sec:stellarator-example}.
Then, in \cref{subsec:Invariance}, we describe our sampling scheme and how we enforce the invariance of the label function. 
In \cref{subsec:RKHS}, we review reproducing kernel Hilbert spaces (RKHS), which define the spaces of functions on which we will approximate our invariant function. 
Each RKHS is equipped with a norm that penalizes non-smooth functions, where the penalization is tunable via the choice of the underlying kernel. 
Finally, we define the invariant boundary value problem (\cref{subsec:IBVP}) and the invariant eigenvalue problem (\cref{subsec:IEP}), which are two natural problems to find the specific smooth nonconstant invariant functions of interest.

\subsection{Preliminaries}
\label{subsec:preliminaries}
We begin by defining some important concepts in this paper. 
We are interested in learning invariant information about a symplectic map $\sympmap : \mathcal X \to \mathcal X$, where $\mathcal X$ is either the plane $\Rbb^2$ or the cylinder $\Tbb \times \Rbb$. 
Here, symplectic means that there is some symplectic 2-form $\omega$ that is invariant under the map, i.e.~$\omega = \sympmap^* \omega$, where $\sympmap^*$ is the pushforward of $\sympmap$. 
Magnetic field-line flow is always symplectic with respect to a form related to conservation of flux.
For the purposes of this paper, the most important quality of symplectic maps is the persistence of invariant circles under perturbation given by the KAM theorem \cite{LLave_2004}. 

Now, let $\lf : \mathcal X \to \Rbb$ be a function of the state space. We say that $\lf$ is invariant under the dynamics if
\begin{equation}
\label{eq:invariance}
    \lf = \mathcal K \lf = \lf \circ \sympmap.
\end{equation}
Here, the linear operator $\mathcal K : \lf \mapsto \lf \circ \sympmap$ is known as the composition or Koopman operator. 
In this way, invariant functions can equivalently be characterized as eigenvectors of the composition operator with eigenvalue $1$. 
In this paper, we will attempt to minimize \eqref{eq:invariance} in a least-squares sense over some region $\Omega$ via an $L^2$ norm
\begin{equation}
\label{eq:approx-invariance}
    \norm{h - h \circ \sympmap}^2_{L^2(\Omega)} = \int_{\Omega} \abs{h(\bm x) - h(\sympmap(\bm x)}^2 \dif \bm x
\end{equation}
Functions for which this quantity is small relative to $H^1$ seminorm $\norm{\nabla h}_{L^2(\Omega}^2$ we will term ``approximately invariant'' or ``label functions.''

One major problem with \cref{eq:approx-invariance} is that constant functions give a zero residual regardless of the map. 
This is the reason for requiring nonconstant label functions, which will be further discussed in the following sections.
When $\sympmap$ is ergodic, the constant function is almost everywhere the unique invariant function.
This fact helps motivate the least-squares approach, as an ergodic map can still have approximately invariant functions.
One consequence is that the existence of an approximately invariant function does not guarantee the existence of a nearby invariant function; such a claim would require a KAM-like result, which is outside the scope of this paper.
In practice, however, we will show there is good agreement between approximately invariant functions and actually invariant structures when such sets exist. 

In contrast, if there are any regions of nested circles within $\Omega$ with positive measure, there are infinitely many linearly independent invariant functions. 
Nested circles are robust to symplectic perturbations in the map by the KAM theorem \cite{LLave_2004}, so we must consider this to be typical. 
Thus, given an infinite number of potential label functions, another major question is what makes for a ``good'' invariant function.
We will opt to find the ``smoothest'' functions that satisfy some boundary condition.

One straightforward method of finding specific invariant functions is via ergodic averaging. 
If we let $f : \mathcal X \to \Rbb$ be an observable (typically in $L^1$), the Birkhoff average $\mathcal B f : \mathcal X \to \Rbb$ is defined as
\begin{equation*}
    \mathcal B f(\bm x) = \lim_{T \to \infty} \frac{1}{T} \sum_{t = 0}^{T-1}(f \circ \sympmap^{t-1})(\bm x).
\end{equation*}
Clearly, if $f$ is invariant then it is a fixed point of the Birkhoff average, i.e.~$\mathcal B f = f$. 
However, it is important to note that the Birkhoff average generically does not satisfy our condition of the invariant function being smooth. 
Even when $f \in C^\infty$ and $\sympmap \in C^\infty$, the Birkhoff average $\mathcal B f$ is typically not even continuous.
For instance, whenever there is a transition from a core to an island chain, one expects to see a jump in the Birkhoff average.

As an algorithm, the naive Birkhoff average is slow to converge, so many iterations of the map $\sympmap$ must be taken. 
If an initial point $\bm x$ is on an invariant circle and $f$ is smooth, the sum can be accelerated using the weighted Birkhoff average \cite{Das2017, sander2020}
\begin{equation}
\label{eq:weighted-birkhoff-average}
    \mathcal{WB} f(\bm x) = \lim_{T \to \infty} \sum_{t = 0}^{T-1} w_{t, T}(f \circ \sympmap^{t-1})(\bm x),
\end{equation}
where $w_{t,T}$ is sampled from a smooth bump function $g : (0,1) \to \Rbb$, such as
\begin{equation*}
    w_{t, T} = \left(\sum_{t = 0}^{T-1} s_{t,T}\right)^{-1} s_{t, T}, \qquad s_{t, T} = g\left(\frac{t+1}{T+1} \right), \qquad g(s) = e^{- \left[s (1-s)\right]^{-1}}.
\end{equation*}
In this situation, the average often converges to machine precision within $10^4$ iterations of the map. 
However, this acceleration does not apply in chaos, as the error empirically goes as $\mathcal O(T^{-1/2})$. 
In \cref{sec:stellarator-example}, the weighted Birkhoff average will be used to verify the label functions that we find.

\subsection{Sampling}
\label{subsec:Invariance}
In this section, we describe the method by which the symplectic map is sampled and how we enforce invariance. 
As an input, we assume that the user provides an oracle for the symplectic map $\sympmap$ and a region of phase space $\Omega$. 
The oracle for $\sympmap$ does not have to be exactly symplectic; each differential equation in this paper is solved by Runge-Kutta integrators, none are exactly symplectic. 
The choice of $\Omega$ is given by the application, but it should fully contain any invariant structure the user wishes to capture. 
For instance, if an invariant circle intersects the boundary of $\Omega$, our method will not have enough information to extrapolate the circle outside the domain. 
In this paper, we always chose a rectangle or an annulus for $\Omega$, but the region could certainly be more complicated.

The method for sampling is simple. 
First, one samples uniformly $N$ points $\bm x_n \in \Omega$. 
Then, $\sympmap$ is applied to each point, resulting in one sequence of $2N$ points $\bm z_n$ where 
\begin{equation}
\label{eq:zn}
    \{ \bm z_n\} = \{\bm x_1, \dots, \bm x_N, \bm F(\bm x_1), \dots, \bm F(\bm x_N)\}.
\end{equation}
For the uniform sampling, we use a low discrepancy Sobol \cite{joe2003} sequence to find $\bm x_n$, but other schemes such as uniformly random sampling and or sampling a grid are also viable. 
We have seen qualitatively worse results in both of these alternative cases, however. 

We choose our sampling scheme because it gives uniform information over a region with a minimal number of applications of $\sympmap$, but it is possible that other sampling methods could be used in other applications. 
For instance, in DMD \cite{kutz2016dynamic}, many short trajectories are often used, rather than single maps.
Similarly, when a region is ergodic, a single long trajectory can be used. 

Now, we consider a candidate label function $\lf$ evaluated on our data. 
The output samples are collected into a vector $\bm \lf$, where $\lf_n = \lf(\bm z_n)$.  
This allows us to discretize \cref{eq:approx-invariance} via the equation
\begin{equation*}
    \Einv = \sum_{n = 1}^{N} \abs{\lf(\sympmap(\bm z_n) - \lf(\bm z_n))}^2 \approx \frac{N}{\mu(\Omega)} \norm{h - h \circ \sympmap}^2_{L^2(\Omega)},
\end{equation*}
where $\mu(\Omega)$ is the Lebesgue measure of the set $\Omega$. We can rewrite this in terms of a matrix equation as
\begin{equation}
\label{eq:Einv}
    \Einv = \norm{\Ginv \bm \lf}^2, \qquad \Ginv = \begin{pmatrix} I & -I \end{pmatrix},
\end{equation}
where $I$ is the $N \times N$ identity matrix. 
We see that $\Ginv$ is rank $N$, so \cref{eq:Einv} has half the number of conditions as the number of unknown values of $\bm \lf$.

\subsection{Reproducing Kernel Hilbert Spaces and Smoothness}
\label{subsec:RKHS}
Now, we turn to the kernel discretization of the label function $\lf$. 
Kernels \cite{schaback2006} are functions on the state space of the form $K : \mathcal X \times \mathcal X \to \Rbb$. 
Often, kernels are chosen to be radial basis functions, meaning $K(\bm x,\bm y)$ only depends on the distance $\norm{\bm x - \bm y}$. 
A popular example is the Gaussian or squared exponential kernel
\begin{equation}
\label{eq:squared-exponential}
    K(\bm x, \bm y) = \exp\left( - \frac{\abs{\bm x - \bm y}^2}{2\sigma^2}\right).
\end{equation}
The squared exponential above, like many kernels, has a tunable width $\sigma$ that defines the natural length scale of $K$.

Given a kernel, the label function can be represented as a linear sum of the kernels centered at the data
\begin{equation}
\label{eq:h-kernel}
    \lf(\bm z) = \sum_{n = 1}^{2N} c_n K(\bm z, \bm z_n).
\end{equation}
Because the representation of $h$ in \cref{eq:h-kernel} is centered at the data, the method in this paper is meshless. 
This means the representation can just as easily be used on arbitrary domains $\Omega$ as rectangular ones.
This is a useful property in the context of stellarators, which often have complicated geometries.
We can evaluate $\lf$ at the nodes, giving
\begin{equation*}
    \bm \lf = \KZZ \bm c, \qquad (\KZZ)_{mn} = K(\bm z_m, \bm z_n).
\end{equation*}
The matrix $\KZZ$ is known as the kernel matrix. 

If $\KZZ$ is always positive definite for distinct sets of points $\bm z_n$, the kernel is called positive definite. 
Positive definite kernels induce a Hilbert space $H_K$, known as the reproducing kernel Hilbert space (RKHS) or native space of $K$. 
Letting $g: \mathcal X \to \Rbb$ be a function of the same form as \eqref{eq:h-kernel} with coefficients $\bm d$, we define the inner product between $\lf$ and $g$ in $H_K$ as
\begin{equation}
\label{eq:RKHS-inner-product}
    \ip{\lf}{g}_K = \bm \lf^T \KZZ^{-1} \bm g = \bm c^T \KZZ \bm d.
\end{equation}
Examining \cref{eq:RKHS-inner-product}, we find evaluation of $\lf$ is equivalent to an inner product $\lf (\bm z_n) = \ip{K(\cdot, \bm z_n)}{h}$. 
This is the ``reproducing'' property for which the RKHS is named.

The norm defined by \eqref{eq:RKHS-inner-product} is 
\begin{equation}
\label{eq:EK}
    \EK = \norm{\lf}_K^2 = \ip{\lf}{\lf}_K = \bm c^T \KZZ \bm c.
\end{equation}
The native RKHS $\mathcal H_K$ of a kernel is given by completing the space of finite kernel representations with respect to $\norm{\cdot}_K$.
The native space of a given kernel is often norm-equivalent to more well known Hilbert spaces such as Sobolev spaces. 
In this way, we use $\norm{\lf}_K$ as a measurement of smoothness, as it will necessarily be large when the derivatives of $h$ are large.
We note that smoothness is an important part of a label function, as the number of invariance conditions is only half the number of the unknowns. 
This means that even in an ergodic region, there are still $N$ functions that exactly satisfy \eqref{eq:Einv}. 
However, these functions tend to overfit the data, and this can be penalized using \eqref{eq:EK}. 
So, the norm in \eqref{eq:EK} can also be seen as penalizing overfitting, and it is commonly used as a regularization term in regression problems. 

The main benefit of kernels for this work is that they offer two major options for tuning: the functional form and the kernel width, also known as the scale factor. 
These options affect both the theoretical qualities of the kernel function and the numerical properties of $\KZZ$.
As the width $\sigma$ increases, the condition number of $\KZZ$ increases, and $\KZZ$ will eventually become numerically low rank.
This effect is more pronounced for smoother kernels such as the squared exponential, whereas sharper kernels such as the inverse multiquadric have slower decay of singular values. 
As such, when the width is higher and the kernel is smoother, the norm in \eqref{eq:EK} will more strongly penalize non-smoothness in $h$. 
We will analyze how these choices affect our method in \cref{sec:stellarator-example}.

\subsection{The Invariant Boundary Value Problem}
\label{subsec:IBVP}
In this section, we consider a problem where we approximately know the locations of an inner and outer invariant circle of a Poincar\'e plot, and we want to know if there are nested invariant circles between those two locations.
In this paper, we represent the domain as an annulus $\Omega = \Tbb \times [a,b]$, where $x_2 = a$ and $x_2 = b$ represent the ``inner'' and ``outer'' circles.
However, any annular domain $\Omega \subset \Rbb^2$ could be used in practice.
We seek a smooth label functions $\lf$ so that we satisfy boundary conditions 
\begin{equation}
\label{eq:strong-boundary-condition}
    \lf(x_1, a) = h_a, \qquad \lf(x_1, b) = h_b,
\end{equation}
where $h_a \neq h_b$.
If $\sympmap$ is smooth and there exists nested invariant circles between the inside and outside the domain, there is an invariant function that can satisfy the boundary conditions. 
Otherwise, if there is a trajectory that approaches both boundaries, it is impossible to exactly find such a function. 
In this way, the boundary value problem does not only ask for an invariant function, but specifically one corresponding to small transport between the inner and outer circles.
We note that while the boundary values $h_a$ and $h_b$ are arbitrary, we require them to be constant on the boundary. 
Otherwise, it may be impossible to find an invariant function that exactly satisfies the boundary conditions. 

Because our sampling scheme does not guarantee points exactly on the boundary $\partial \Omega$, we require a more general notion of boundary conditions.
So, instead of strictly enforcing boundary conditions as constraints, we consider defining two functions on our state space: a boundary value function $\hbd: \Tbb \times \Rbb \to \Rbb$ and a weighting function $\wbd : \Tbb \times \Rbb \to \Rbb^{\geq 0}$.
The function $\hbd$ encodes the value of the boundary condition, and $\wbd$ indicates how strongly the boundary condition should be enforced.
Given the boundary functions, we define a boundary energy
\begin{equation}
\label{eq:Ebd}
    \Ebd = \norm{\bm \lf - \bm \hbd}_{\Wbd}^2 = (\bm \lf - \bm \hbd)^T \Wbd (\bm \lf - \bm \hbd),
\end{equation}
where $(\bm \hbd)_n = \hbd(\bm z_n)$ and $\Wbd$ is a diagonal matrix with $(\Wbd)_{nn} = \wbd(\bm z_n)$.
Clearly, \cref{eq:Ebd} is zero if the boundary conditions are exactly satisfied on $\bm z_n$, and positive when they are not. 

A simple example of functions $\hbd$ and $\wbd$ that weakly enforce \cref{eq:strong-boundary-condition} is to define two buffer regions where the boundary conditions are enforced: $\Gamma_a = \{(x,y) \in \Tbb \times \Rbb \mid y < a + \beta\}$ and $\Gamma_b = \{(x,y) \in \Tbb \times \Rbb \mid y > b - \beta\}$ for some buffer $\beta \geq 0$. 
Then, the boundary functions are
\begin{equation*}
    \hbd = h_a \mathbbm{1}_{\Gamma_a} + h_b \mathbbm{1}_{\Gamma_b}, \qquad \wbd = \mathbbm{1}_{\Gamma_a} + \mathbbm{1}_{\Gamma_b},
\end{equation*}
where $\mathbbm{1}_{\Gamma}$ is the indicator function 
\begin{equation*}
    \mathbbm{1}_\Gamma(\bm x) = 
    \begin{cases} 
        1, & \bm x \in \Gamma, \\ 
        0, & \bm x \notin \Gamma. 
    \end{cases}
\end{equation*}
These functions essentially force a strip of the label function to take the boundary value near the boundary.
The boundary conditions also apply to points that may have been mapped out of $\Omega$, addressing potential problems of not requiring $\Omega$ to be forward-invariant. 

In practice, we have found that smoothing $\wbd$ near the boundaries of $\Gamma_a$ and $\Gamma_b$ is often useful. 
This allows us to more weakly enforce the boundary condition when points are close to the edge, as well as making the energy smooth under perturbations in $\bm z_n$. 
Consequently, if $\sympmap$ depends smoothly on a parameter $p$, then $\Ebd$ is also smooth in $p$ (see \cref{sec:standard-map}). 
As an example, consider the boundary functions
\begin{align}
\label{eq:hbd}
    \hbd &= \frac{h_a + h_b}{2} + \frac{h_b - h_a}{2} \tanh\left( \frac{2y - a - b}{2\alpha}\right), \\
\label{eq:wbd}
    \wbd &= \Sigma\left(\frac{y - b + \beta}{\alpha}\right) + \Sigma\left(- \frac{y - a - \beta}{\alpha}\right),
\end{align}
where $\Sigma$ is the sigmoid function $\Sigma(y) = \frac{1}{1 + e^{-y}}$, $\beta$ is a buffer length, and $\alpha > 0$ is a smoothing width.
Equations \cref{eq:hbd} and \eqref{eq:wbd} can be thought of as smoothing the sharp $\Gamma$ boundary conditions.

We are now ready to define our boundary value problem. We attempt to minimize the following residual:
\begin{align}
\nonumber
    R &= \min_{h\in H_K} \left(\Ebd + \Einv + \epsilon \EK\right), \\
\label{eq:bvp-orig}
      &= \min_{h \in H_K} \left(\norm{\bm \lf - \bm \hbd}_{\Wbd}^2 + \norm{\Ginv \bm h}^2 + \epsilon \norm{h}_K^2\right),
\end{align}
where $0 < \epsilon \ll 1$ is a small parameter and $\Ebd$, $\Einv$, and $\EK$ are defined in \cref{eq:Ebd}, \cref{eq:Einv}, and \cref{eq:EK} respectively. 
The smallness of $\epsilon$ means that the invariance and boundary energies are treated as weakly enforced constraints with regularization. 
The residual $R$ is a measure of ``goodness of fit.'' 
Assuming the kernel is wide enough, small values of $R$ indicate that we have found a good label function.

Equation \cref{eq:bvp-orig} satisfies the conditions of the representer theorem in \cite{scholkopf2001}, meaning the solution $h$ must take the form \cref{eq:h-kernel}. 
Thus, we can rewrite \eqref{eq:bvp-orig} as a finite dimensional minimization problem
\begin{equation}
\label{eq:bvp}
    R = \min_{\bm c \in \Rbb^{2N}} \left(\norm{\KZZ \bm c - \bm \hbd }_{\Wbd}^2 + \norm{\Ginv \KZZ \bm c}^2 + \epsilon \bm c^T \KZZ \bm c\right),
\end{equation}
We can solve \cref{eq:bvp} by setting the derivative with respect to $\bm c$ equal to zero and left-multiplying by $\KZZ^{-1}$ to find
\begin{equation}
\label{eq:bvp-solve}
    \left((\Wbd + \Ginv^T \Ginv) \KZZ + \epsilon I \right) \bm c = \Wbd \bm \hbd.
\end{equation}
We can solve this equation via dense linear solves, leading to an asymptotic cost of $\mathcal O(N^3)$. 
In summary, the process for finding an approximately invariant label function is given by the following pseudocode:
\begin{algorithm}[H]
\caption{Invariant Boundary Value Problem}
\label{alg:IBVP}
\begin{algorithmic}[1]
\Require Symplectic map $\sympmap$, domain $\Omega$, kernel $K$, number of points $N$, regularization $\epsilon$, and boundary functions $\hbd$ and $\wbd$
\State Sample $\bm z_n \in \Omega$ uniformly for $1 \leq n \leq N$
\State $\bm z_{N+n} \gets \bm \sympmap(\bm z_n)$ for $1 \leq n \leq N$
\State Sample boundary functions $\hbd$ and $\wbd$ at each $\bm z_n$ and form $\KZZ$
\State Solve \cref{eq:bvp-solve} for $\bm c$ and compute $R$
\Ensure Approximate label function $h$ and residual $R$
\end{algorithmic}
\end{algorithm}

As an example, consider the nonlinear pendulum on $\Tbb \times \Rbb$ with the dynamics
\begin{equation}
\label{eq:pendulumODE}
    \dot x = y, \qquad \dot y = - \sin (2\pi x).
\end{equation}
We obtain a symplectic map $\sympmap$ by evolving the above differential equations a time period of $\sqrt 2$. 
For $\bm z_n$, uniformly sample the domain $\Omega = \Tbb \times [-2.1, 2.1]$. 
For boundary functions, we use \cref{eq:hbd} and \cref{eq:wbd} with $\alpha = 0.02$, $\beta = 0.1$, $h_a = -1$ and $h_b = 1$. 
Because the dynamics are periodic in $x$, we use a product kernel with a squared sine exponential and a squared exponential, i.e. 
\begin{equation}
\label{eq:sine-kernel}
    K(\bm x, \bm x') = \exp\left(-\frac{\sin^2(\pi \abs{x - x'})}{2 \pi \sigma^2} - \frac{\abs{y - y'}^2}{2\sigma^2}\right),
\end{equation}
where the kernel width $\sigma$ takes a value of $0.5$. 
We sample $N = 100$ points and regularize with $\epsilon = 10^{-8}$.

\begin{figure}
    \centering
    \includegraphics[width=1.0\textwidth]{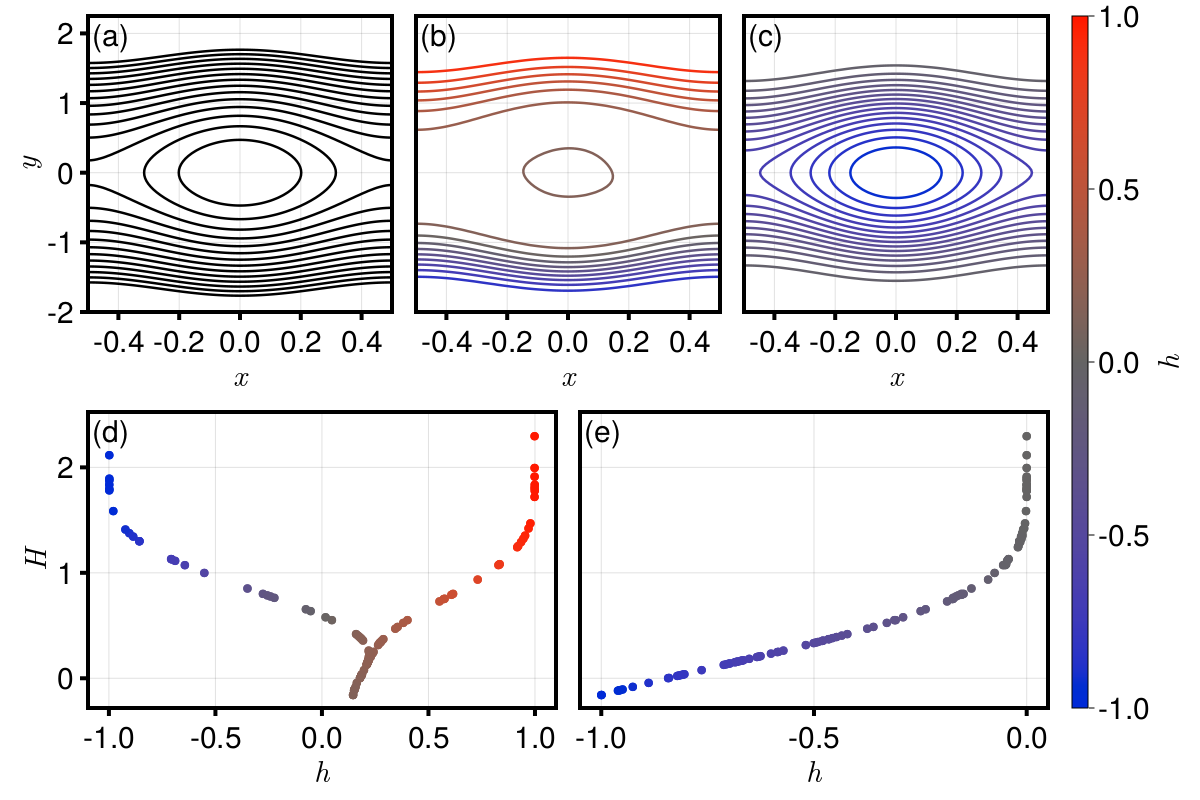}
    \caption{(a) Level sets of the Hamiltonian for the nonlinear pendulum \cref{eq:H_pendulum}, (b) a label function obtained by solving the boundary value problem \eqref{eq:bvp} with $N=100$ evaluations of the associated symplectic map, and (c) a label function obtained by solving the eigenvalue problem \eqref{eq:rayleighquotient} with the same parameters. The contour shapes of all three plots are qualitatively similar. We plot the Hamiltonian $H(\bm x_n)$ vs (d) the boundary value problem and (e) the eigenvalue problem label functions $h(\bm x_n)$ on the sampled data. We see the Hamiltonian is locally a plot of both functions, indicating the label functions are learning the correct invariant sets. }
    \label{fig:nonlinear-pendulum}
\end{figure}

The solution to this problem is shown in \cref{fig:nonlinear-pendulum} (b). We see that there is a large gradient in $h$ where there are nested flux surfaces, whereas $h$ is approximately flat in the center ``island.''
This is because while the boundary conditions force the label function to have a gradient from bottom to top, no progress to this end can be made while staying invariant on the separatrix. 
Additionally, the smoothness energy otherwise encourages flatness, so the island converges to a flat value. 
The residual of the least-squares problem is $R = 2.57*10^{-7}$, with $\EK = 18.2$, $\Ebd = 2.58 \times 10^{-7}$ and $\Einv = 6.54 \times 10^{-9}$. 
The value of $\Einv$ tells us that the function is very close to invariant at the sampled data. 

In \cref{fig:nonlinear-pendulum} (a), we plot the level sets of the Hamiltonian
\begin{equation}
\label{eq:H_pendulum}
    H = \frac{1}{2} y^2 - \frac{1}{2\pi} \cos(2 \pi x)
\end{equation}
for comparison. 
Because the Hamiltonian is invariant under the flow, we expect that the contours of the match the contours of the label function.
Visually, we see good agreement where there are nested circles, but they are farther from correct in the center where $h$ is approximately flat.
Moreover, in regions of nested circles where $h$ has a nonzero derivative, we expect $H(\bm x)$ to be approximately a function of $h(\bm x)$.
We plot $h(\bm x_n)$ vs $H(\bm x_n)$ in \cref{fig:nonlinear-pendulum} (d), showing a functional relationship between the values in regions of nested circles (i.e.~where $-1 < h < 0.1$ and $0.3 < h < 1$). 
In the regions where $h$ is about flat (i.e.~the boundaries $h \approx \pm 1$ and in the island $h \approx 0.2$), the function becomes close to vertical. 
We note the $H$ is not symmetric to sign changes in $h$ due to the asymmetry of the sampling of $\bm x_n$ in the domain.
The boundary value problem will be further explored via the standard map in \cref{sec:standard-map}.  

\subsection{The Invariant Eigenvalue Problem}
\label{subsec:IEP}
We now describe an alternative approach to the previous section. 
We consider the case where we only know the boundary of a stellarator configuration and we want to find invariant functions inside.
Because we do not have specific information about the inside of the domain, we do not have a natural second boundary through which we can force the label function to be nonzero. 
For this reason, we instead consider the Rayleigh quotient
\begin{equation}
\label{eq:rayleighquotient}
    \lambda = \min_{\bm c \neq 0} \frac{\norm{\Ginv \KZZ \bm c}^2 + \norm{\KZZ \bm c}^2_{\Wbd} + \epsilon \norm{\bm c}^2_K}{\norm{\KZZ \bm c}^2}.
\end{equation}
We recognize the numerator terms, from left to right, as the invariant energy \cref{eq:Einv}, a zero boundary condition \cref{eq:Ebd}, and the regularizing kernel norm \cref{eq:EK} where $0 < \epsilon \ll 1$. 
In the denominator, we have made the problem scale invariant by normalized by the $L^2$ norm of the label function, allowing for a nonzero solution to the problem.
Additionally, due to the boundary condition penalty, we remove the constant function as a solution.

The Rayleigh quotient \cref{eq:rayleighquotient} is minimized by the smallest eigenvalue $\lambda$ of the eigenvalue problem
\begin{equation}
\label{eq:eval-bad}
    \left( \KZZ ( \Ginv^T \Ginv + \Wbd) \KZZ + \epsilon \KZZ \right) \bm c = \lambda \KZZ^2 \bm c.
\end{equation}
However, $\KZZ$ is numerically low rank for large $\sigma$, so \cref{eq:eval-bad} is difficult to solve directly.
To isolate the dependence on $\KZZ$, we introduce the auxiliary variable $\bm h = \KZZ \bm c$, giving
\begin{equation*}
    \begin{pmatrix}
        \Ginv^T \Ginv + \Wbd + \delta I & I \\
        I & -\epsilon^{-1} \KZZ
    \end{pmatrix}\begin{pmatrix}
        \bm h \\ \epsilon \bm c
    \end{pmatrix} = (\lambda + \delta) \begin{pmatrix}
        \bm h \\ 0
    \end{pmatrix},
\end{equation*}
where we have shifted the eigenvalue problem by a fixed amount $0 < \delta \ll 1$. To avoid direct inversions of $\KZZ$, we invert the above system by a Schur complement, giving
\begin{align}
\label{eq:eval}
     (\lambda + \delta)^{-1} \bm h &= \left(A^{-1} - A^{-1}(\epsilon^{-1} \KZZ + A^{-1})^{-1}A^{-1} \right)\bm h, \\
\nonumber
     \epsilon \bm c &=  (\epsilon^{-1} \KZZ + A^{-1})^{-1}A^{-1} \bm h,
\end{align}
where $A = \Ginv^T \Ginv + \Wbd + \delta I$. We see that $A^{-1}$ regularizes the inversion of $\KZZ$, reducing ill-conditioning of the problem. 

To solve the shifted and inverted eigenvalue problem \cref{eq:eval}, we must invert multiple operators. The matrix $A$ is block $2\times 2$ diagonal, so the computation of $A^{-1}$ is performed directly in $\mathcal{O}(N)$ time. Then, we precompute the Cholesky factorization $U^T U = \epsilon^{-1} \KZZ + A^{-1}$. At $\mathcal{O}(N^3)$ time, this is the dominant cost of the algorithm. Finally, we solve the eigenvalue problem for the largest eigenvalue $(\lambda + \delta)^{-1}$ using an Arnoldi method implemented in the Julia wrapper of ARPACK.
The value of $\bm c$ is obtained using the same precomputed matrices.
The pseudocode for this algorithm is the following:
\begin{algorithm}[H]
\caption{Invariant Eigenvalue Problem}
\label{alg:IEVP}
\begin{algorithmic}[1]
\Require Symplectic map $\sympmap$, domain $\Omega$, kernel $K$, number of points $N$, regularization $\epsilon$, boundary functions $\wbd$, shift parameter $\delta$
\State Sample $\bm z_n \in \Omega$ uniformly for $1 \leq n \leq N$
\State $\bm z_{N+n} \gets \bm \sympmap(\bm z_n)$ for $1 \leq n \leq N$
\State Sample boundary functions $\hbd$ and $\wbd$ at each $\bm z_n$ and form $\KZZ$
\State Directly compute $A^{-1} = (\Ginv^T \Ginv + \Wbd + \delta I)^{-1}$
\State Cholesky factorize $U^T U \gets \epsilon^{-1} \KZZ + A^{-1}$ 
\State Solve \cref{eq:eval} for the largest eigenvalue $(\lambda+\delta)^{-1}$ and eigenvectors
\Ensure Smallest eigenvalue $\lambda$ and label eigenfunction $\bm h$
\end{algorithmic}
\end{algorithm}

We show an example of the eigenvalue problem \cref{eq:rayleighquotient} applied to the pendulum example \cref{eq:pendulumODE} in \cref{fig:nonlinear-pendulum} (c) using the same parameters as \cref{subsec:IBVP}. The resulting eigenvalue is $\lambda = 3.905\times 10^{-10}$, showing the eigenfunction strongly satisfies the boundary conditions and invariance at the sampled points. Again, we see that the contour shapes qualitatively match both the Hamiltonian (\cref{fig:nonlinear-pendulum} (a)) and the boundary value problem version of the problem (\cref{fig:nonlinear-pendulum} (c)). However, the recovered eigenfunction emphasizes the center much more than the boundary value problem. This is because the center looks akin to a ``core'' to the eigenvalue problem, and allows for more concentration of the mass of the eigenfunction.
In \cref{fig:nonlinear-pendulum} (e), we plot the eigenfunction $h(\bm x_n)$ against the Hamiltonian $H(\bm x_n)$, showing that $H$ can be found approximately as a function of the label function.

\section{Example: A Stellarator Map}
\label{sec:stellarator-example}

\begin{figure}
    \centering
    \includegraphics[width=0.75\textwidth]{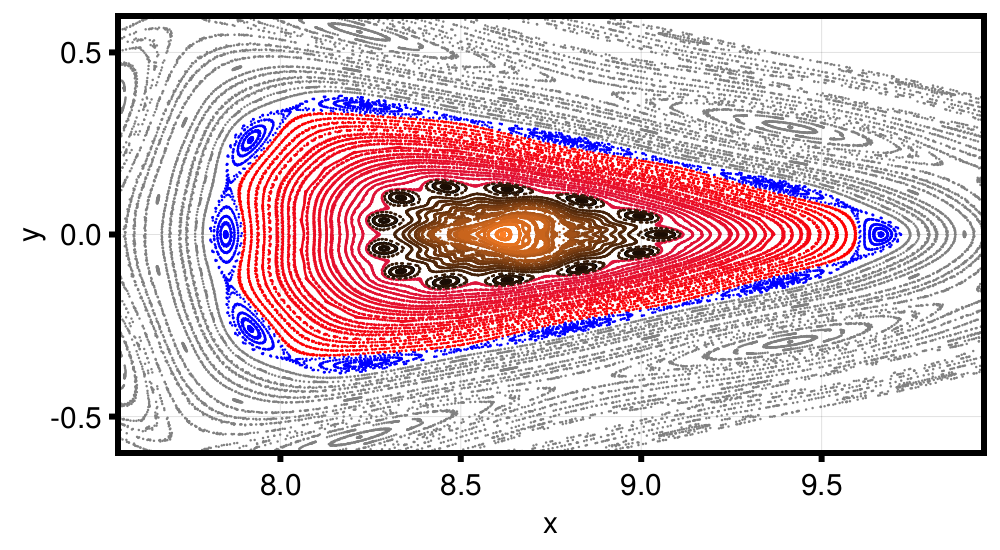}
    \caption{
    A Poincar\'e plot of a candidate stellarator configuration. The plot is colored into four regions: (blue) a 10-island chain, (red) a region of nested invariant circles, (brown) a core consisting of islands, nested circles, and small amounts of chaos, and (gray) an outer region. Compare with the label functions in \cref{fig:stellarator_label_functions}. The Poincar\'e plot consists of $80$ trajectories with $1001$ points each, resulting in $80000$ evaluations of the symplectic map over 81.5 seconds.
    }
    \label{fig:stellarator_poincare_plot}
\end{figure}

In this section, we analyze the Rayleigh quotient method \cref{eq:rayleighquotient} using an example stellarator return map, plotted in \cref{fig:stellarator_poincare_plot}. 
The symplectic map is obtained by evolving the ODE 
\begin{equation}
\label{eq:B-field-tracing}
    \dot {\bm x} = \frac{1}{B^\phi(\phi, \bm x)}\bm B(\phi, \bm x), \qquad 0 \leq \phi \leq \frac{2 \pi}{n_{\mathrm{fp}}},
\end{equation} 
where $\bm x(\phi) = (x(\phi),y(\phi))$ is the 2-vector of radial and vertical positions in a cross-section with fixed poloidal angle $\phi$, $\bm B = (B^x, B^y)$ is in-plane magnetic field, $B^\phi$ is the out-of-plane magnetic field, and $n_{\mathrm{fp}}$ is the number of field periods of the stellarator configuration. 
The ODE \cref{eq:B-field-tracing} is solved by a Runge-Kutta 4(5) integrator.
This particular stellarator was chosen primarily for its combination of islands on the edge and its complicated core, consisting of islands, nested circles, and chaos.
In between the islands and the core, there is a large region of nested flux surfaces.
Due to all of these features, the map is more complicated than the perturbed and non-perturbed pendulums in the previous section, giving a more realistic test case.
This example also gives useful timing results.
As a baseline for the timing, the plot in \cref{fig:stellarator_poincare_plot} required $80000$ evaluations of the map over $81.5$ seconds, resulting in an average symplectic map evaluation taking around $1$ millisecond. 

For this example, we are primarily interested in exploring the dependence of the label function invariance on hyperparameters, including the number of evaluations $N$, the kernel type $K$ and width $\sigma$, and the regularization parameter $\epsilon$. 
As such, we need a method for validating the invariance of a given map. 
We note that while $\Einv$ may be small for a given label function, it is possible that the map is not actually invariant off the sampled data. 
Additionally, $\Einv$ is only an indicator of how invariant the map is over a single iteration, so it is not necessarily a good indicator of how invariant a trajectory is. 

As an alternative, we consider taking a weighted Birkhoff average \cref{eq:weighted-birkhoff-average} of our label function $h$ as a validation error. 
As we mentioned in \cref{subsec:preliminaries}, if $h$ is invariant under the flow, then $\WB{h} = h$.
Otherwise, the values will deviate, and we consider the label function to be good if $\WB{h} \approx h$.
As a measurement of this, we use the normalized sum-of-squares residual
\begin{equation}
 \label{eq:S}
    S = \frac{1}{S_0} \sum_{j = 1}^J (h(\bm x_j) - \mathcal{WB}[h](\bm x_j))^2, \qquad S_0 =  \sum_{j = 1}^J (h(\bm x_j) - \bar h)^2, 
\end{equation}
where the $J$ points $\bm x_j$ are uniformly sampled over the domain and $\bar h$ is the mean of the sampled label function values $h(\bm x_j)$. 
The sum-of-squares residual is related to the $R^2$ coefficient of determination of the line $\WB{h} = h$ via $S = 1 - R^2$. 
Throughout this section, we use $J = 1000$ points and a Birkhoff average length of $T = 100$ points. 
Due to the large number of $\sympmap$ and $h$ evaluations required to compute the Birkhoff average, this process is likely too slow for practical situations, and is primarily useful as a diagnostic tool. 

\begin{figure}
    \centering
    \includegraphics[width=0.95\textwidth]{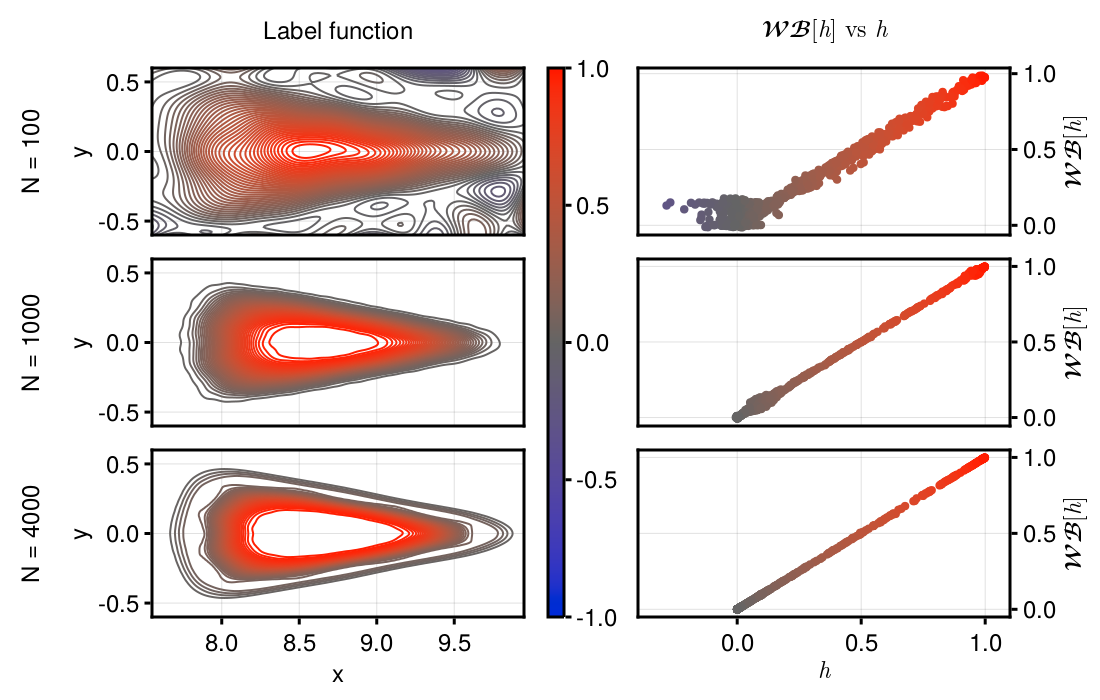}
    \caption{
        (left) Three label functions $h$ for the example stellarator map with $N = 100$, $1000$, and $4000$, $\sigma$ scaling as $2.83 / \sqrt{N}$, $K$ as the squared exponential \cref{eq:squared-exponential}, and $\epsilon = 10^{-8}$. The functions $h$ are normalized to have a maximum absolute value of $1$. We see that the $N=100$ plot captures the general shape of the label function with large error near the boundary, $N=1000$ starts to resolve the chaos in the center, and $N=4000$ resolves the finer island structures. (right) Plots of the weighted Birkhoff average $\WB{h}$ averaged over $T=100$ points vs the label function $h$ for $J=1000$ initial points $\bm x_j$. The points are colored by $h$ to match the left column. The correlation between $h$ and $\WB{h}$ improves with $N$, with a validation error \cref{eq:S} of $S = 7.37 \times 10^{-2}$ for $N=100$, $S = 6.11 \times 10^{-4}$ for $N=1000$, and $S = 9.41 \times 10^{-6}$ for $N=4000$. The unresolved islands in the $N=1000$ label function are represented by two bulges in $\WB{h}$ vs $h$ near $h = 0.1$ and $h = 0.95$. 
    }
    \label{fig:stellarator_label_functions}
\end{figure}

As an example of the sum-of-squares residual, we consider the three label functions depicted in \cref{fig:stellarator_label_functions}. 
We find the label functions by solving the eigenvalue problem \cref{eq:rayleighquotient} on the domain $\Omega = [7.5, 10] \times [-0.65, 0.65]$ with a boundary condition $\wbd = \mathbbm{1}_\Gamma$ where $\Gamma = \Rbb^2 \backslash \left( [7.55, 9.95] \times [-0.6, 0.6]\right)$. 
We use the same domain and boundary conditions throughout this section. 
Additionally, we consider $\epsilon = 10^{-8}$, three values of $N\in \{100, 1000, 4000\}$, and $\sigma$ scaling with $N$ via the density of points as 
 \begin{equation}
 \label{eq:sigma-scaling}
     \sigma = \frac{\sigma_0}{\sqrt N},
 \end{equation}
where $\sigma_0 = 2.83$ in this example.

The resemblance of the label functions in \cref{fig:stellarator_label_functions} compared to the Poincar\'e plot in \cref{fig:stellarator_poincare_plot} improves as $N$ increases. 
For $N = 100$ (\cref{fig:stellarator_label_functions} top left), the shape is correct, but the label function does not resolve chaotic core or the outer island ring. Additionally, there are not enough points for this kernel to resolve the boundary region to zero.
In the right column we plot $\WB{h}$ vs $h$, the inputs of the validation metric \cref{eq:S}.
The mismatches between the Poincar\'e plot and the label function appear in the \cref{fig:stellarator_label_functions} (upper right), where it is clear that $\WB{h} \approx h$, but the line is blurry. 
The resulting value of $S$ is $7.37 \times 10^{-2}$.

For $N = 1000$ (\cref{fig:stellarator_label_functions} middle left), the label function matches the Poincar\'e plot in finer detail by resolving the boundary and flattening the core. 
However, it still does not resolve the outer or inner rings of islands, resulting in two bulges in the $\WB{h}$ vs $h$ plot and a validation error of $S = 6.11 \times 10^{-4}$. 
For $N = 4000$, the island rings are resolved, the plot of $\WB{h}$ vs $h$ is visually thin and straight, and the value of $S$ decreases further to $9.41 \times 10^{-6}$.

\begin{figure}
    \centering
    \includegraphics[width=0.9\textwidth]{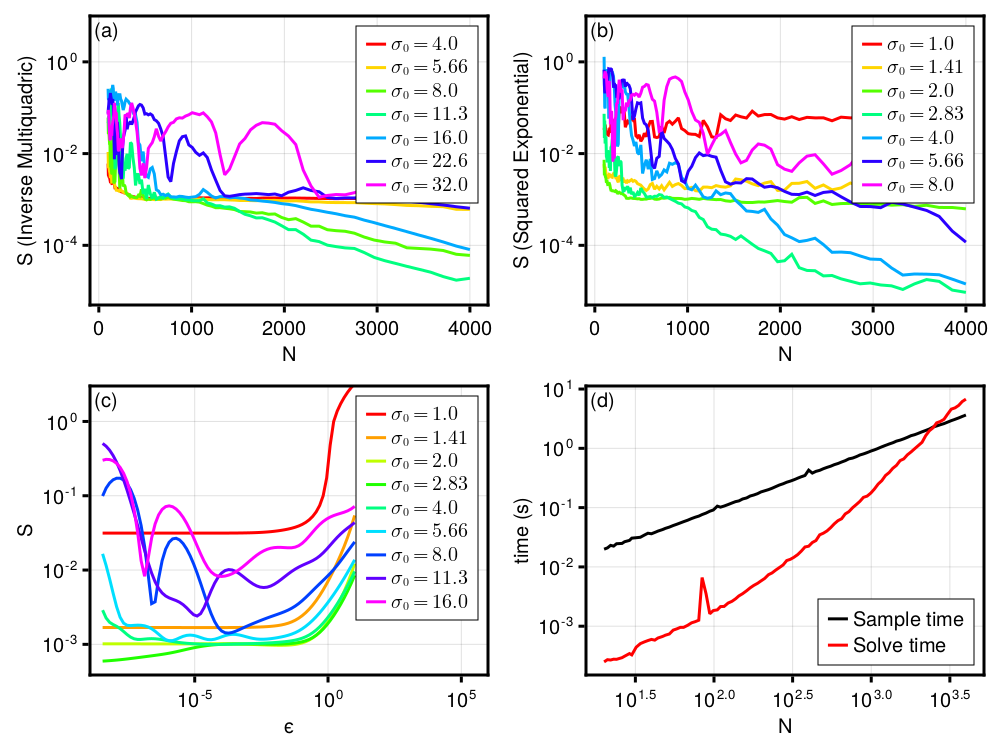}
    \caption{
        (a,b) Plots of the validation error \cref{eq:S} vs $\sigma$ \cref{eq:sigma-scaling} for varying $\sigma_0$. The kernels used are (a) the inverse multiquadric \cref{eq:inversemultiquadric} and (b) the squared exponential \cref{eq:squared-exponential}. We see the two kernels achieve similar errors for equal values of $N$. (c) A plot showing the validation error dependence on $\epsilon$ for the squared exponential kernel with $N=100$ and varying $\sigma_0$. (d) A timing comparison between the time of sampling and the time of solving for the kernel function as a function of $N$ for the squared exponential kernel with $\sigma_0 = 2.83$ for a single core. The two times match order at around $N=2500$. 
    }
    \label{fig:kernel_width}
\end{figure}

Now, we explore the relationship between $N$, $\sigma$, $K$, and $S$ in more detail. 
We choose two kernels for this experiment: the squared exponential \cref{eq:squared-exponential} as an example of a smooth localized kernel and the inverse multiquadric 
\begin{equation}
\label{eq:inversemultiquadric}
    K(\bm x, \bm y) = \frac{1}{\sqrt{1 + |\bm x - \bm y|^2/\sigma^2}}
\end{equation}
as an example of a sharp heavy-tailed kernel.
For each kernel, we choose $7$ values of $\sigma_0$ and $100$ values of $N$ between $100$ and $4000$, both logarithmically distributed. 
In every trial, we have $\epsilon = 10^{-8}$ fixed. 
Then, for every combination of $N$, $\sigma_0$, and $K$, we compute the kernel function via \cref{eq:rayleighquotient} and $S$.

\cref{fig:kernel_width} (a) and (b) show the results for the inverse multiquadric and squared exponential respectively. 
For both kernels, there are specific values of $\sigma_0$ that are best for large $N$, with $\sigma_0 = 11.3$ for the inverse multiquadric and $\sigma_0 = 2.83$ for the squared exponential. 
The best validation error for the two kernels is approximately equal for fixed $N$, giving little indication that one kernel is better than another in this application.
For both kernels, when $\sigma_0$ is chosen too large the label functions converge, but much slower than the optimal value. 
When $\sigma_0$ is too small, both kernels fail to converge, as the label function tends to overfit to the invariance energy.

In \cref{fig:kernel_width} (c), we consider the dependence of $S$ upon $\epsilon$ for the squared exponential kernel. 
We again scan over values of $\sigma_0$, except now we fix $N = 1000$ and scan $70$ values of $\epsilon$ logarithmically spaced between $10^{-8}$ and $10$.
The dependence on $\epsilon$ takes two forms depending on the value of $\sigma_0$. 
For the values of $\sigma_0 \geq 4.0$, there is a clear minimum in the value of $S$ vs $\epsilon$. 
This is due to the kernel matrix being low rank for wide kernel widths, forcing the regularization to take a larger role.
For $\sigma_0 \leq 2.0$, the dependence of $S$ on $\epsilon$ is flat for very large regions, showing a lack of sensitivity to the parameter. \cref{fig:kernel_width} (b) of $\sigma = 2.83$, the value of $S$ actually improves when $\epsilon$ is decreased from $10^{-5}$ to $10^{-9}$, showing that a well tuned value of $\sigma$ and very small values of $\epsilon$ are required for accurate label functions.

Finally, in \cref{fig:kernel_width} (d), we consider the timing of the method. 
The code was timed on a personal computer on a single core with an Intel Core i7-7700 3.6 GHz processor and 16 GB RAM. 
We time both the sampling phase and the solve phase of \cref{alg:IEVP} for the $\sigma_0 = 2.83$ line of \cref{fig:kernel_width} (b). 
The leading order cost of the solve phase is $\mathcal O(N^3)$, due to the Cholesky factorization and eigenvalue solve.
Until around $N=1000$ the $\mathcal O(N)$ cost of sampling outweighs the cost of the label function retrieval, each of which takes around a second at $N=1000$.
Beyond that point, the total time to compute the accurate $N=4000$ label function is $10.4$ seconds combined, which is still about an eighth of the time to compute the Poincar\'e plot.

\section{Example: The Standard Map}
\label{sec:standard-map}

\begin{figure}
    \centering
    \includegraphics[width=0.9\textwidth]{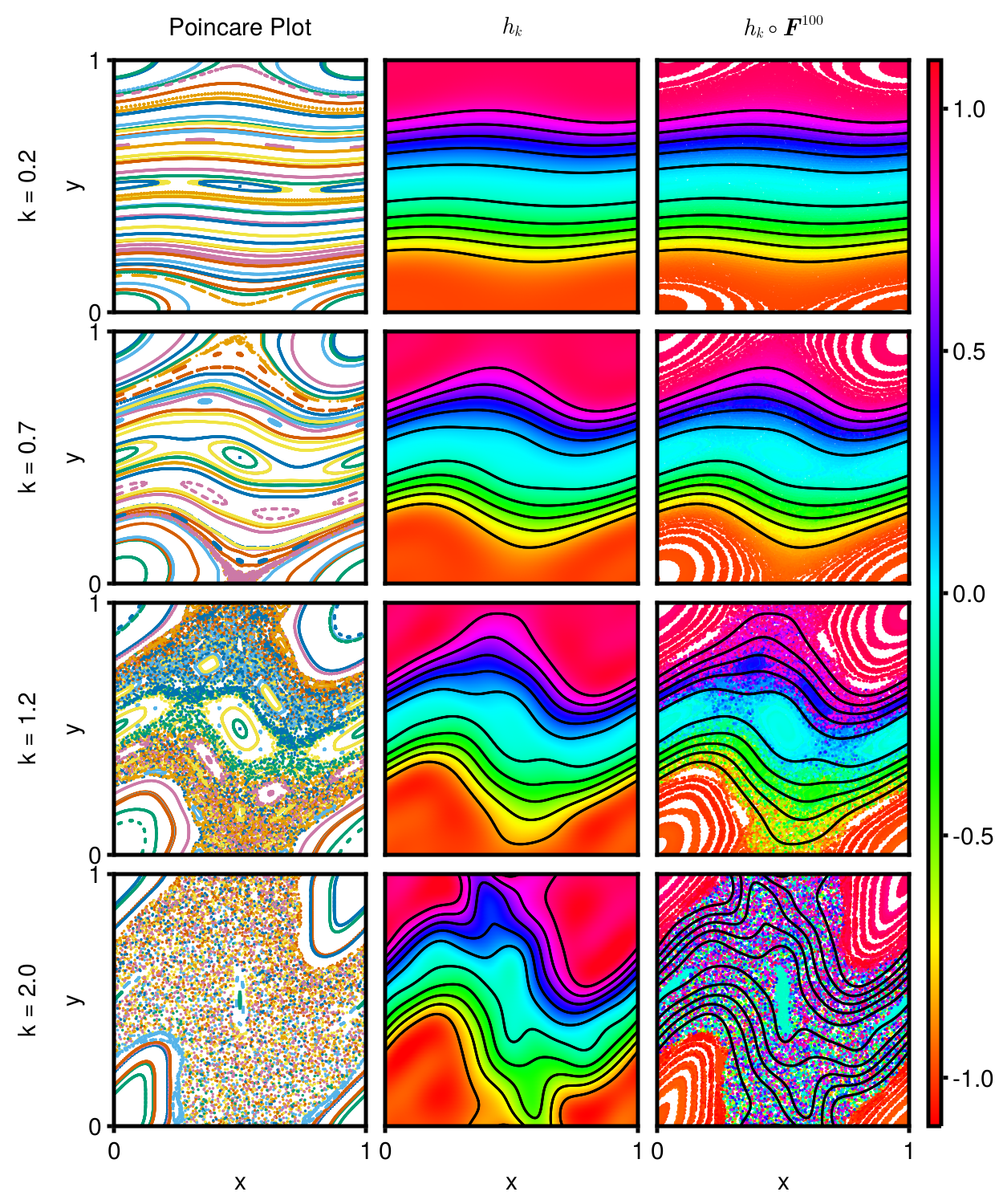}
    \caption{
    (left) Poincar\'e plots of the standard map \cref{eq:standard-map} for $k = 0.2$, $0.7$, $1.2$, and $2.0$. As $k$ increases, the map has more islands and chaos.
    (middle) A label function computed for the standard map with a product squared exponential kernel \eqref{eq:sine-kernel}, $\sigma = 0.2$, $N=500$, and $\epsilon = 10^{-5}$. The black contours are level sets of the label function.
    (right) The middle label function composed with $100$ iterations of the standard map. For $k=0.2$ and $0.7$, the label function is approximately equal, with the only errors due to island mixing. At $k=1.2$, the chaos causes slow diffusion across field lines. For $k=2.0$, the chaos fully mixes the points.
    }
    \label{fig:chaos-functions}
\end{figure}

In this section, we give an example of the boundary value method \cref{eq:bvp} for the standard map. The standard map is defined via the iteration
\begin{align}
\nonumber
    a_{t+1} &= a_t + b_{t+1},\\
\label{eq:standard-map}
    b_{t+1} &= b_t - \frac{k}{2\pi} \sin (2\pi a_t),
\end{align}
where $k$ is a free parameter. 
In this way, we define the symplectic map on $\Tbb \times \Rbb$ as $\sympmap_k : (a_t, b_t) \mapsto (a_{t+1}, b_{t+1})$. 
At $k = 0$, the dynamical system is integrable with $b_t = b_0$ as a constant of motion, meaning the entire state space is foliated by invariant circles. 
When $k$ increases, the invariant circles break down, until eventually the final circle disappears at a critical $k_c \approx 0.971635$ \cite{greene1979,haro2016}, allowing for transport throughout the state space. 
Once $k$ surpasses $k_c$ trajectories can pass from $y = 0$ to $y = 1$, but there are still many islands in state space. 
As $k$ continues to increase, these islands diminish as the volume and rates of chaos increase.
One of the main ways of studying transport in the standard map is via turnstiles and cantori \cite{meiss2015}.
Poincar\'e plots for $k = 0.2$, $0.7$, $1.2$, and $2.0$ are all shown in the left column of \cref{fig:chaos-functions}. 

The changes in the standard map with $k$ give a natural way to analyze at the dependence of \cref{alg:IBVP} on the map $\sympmap_k$.
This is in contrast to the previous section, where we were primarily interested in the dependence on hyperparameters. 
To look at dependence on $k$, we consider solving \cref{eq:bvp} on the domain $\Omega = \Tbb \times [0,1]$ with $N = 500$ and the product kernel \cref{eq:sine-kernel} with $\sigma = 0.2$. 
Throughout this section, we use the regularization $\epsilon = 10^{-5}$ and the smoothed boundary conditions \cref{eq:hbd,eq:wbd} with $\alpha = \beta = 1/100$, $\lf_a = -1$, and $\lf_b = 1$.

We plot four label functions $\lf_k$ for $k\in \{0.2, 0.7, 1.2, 2.0\}$ in the middle column of \cref{fig:chaos-functions}. 
The plots are obtained by first scattering a grid of points on $\Omega$ colored by the value of $\lf_k$ at those points.
On top of the grid, we plot contours of $\lf_k$ in black.
In the right column of \cref{fig:chaos-functions}, we advect the points in the middle column via $100$ iterations of $\sympmap_k$, 
keeping the original contours.
This gives a visual of how the label function is mixed under the map.

For $k = 0.2$ and $0.7$, the label functions $\lf_k$ and $\lf_k \circ \sympmap_k^{100}$ both qualitatively match the Poincar\'e plot in the left column. 
In both cases, the label function has gradients in regions with nested circles and is more flat in regions with islands. 
The flat regions make it so that the advected color is difficult to distinguish from the original plot.
For $k = 1.2$, the label function still appears to match the dynamics shown in the Poincar\'e plot, despite the obviously chaotic trajectories.
In the advected plot, we see that while the trajectories cross the black contours of $h$, the diffusion of points is slow, with the advected plot retaining the impression of the original.
In the final $k = 2.0$ plot, the obtained label function has little resemblance to the Poincar\'e plot, and this is confirmed by the advected plot on the right.

\begin{figure}
    \centering
    \includegraphics[width=0.95\textwidth]{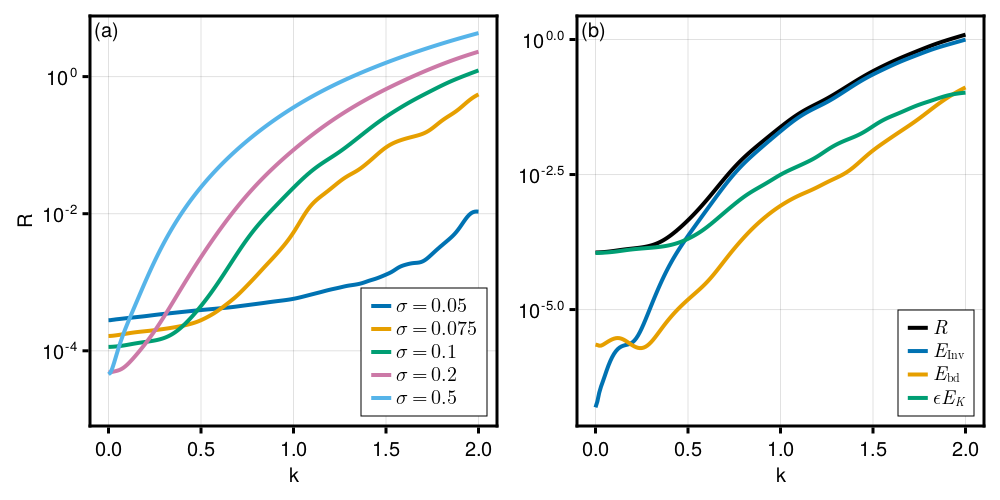}
    \caption{
    (a) The residual \cref{eq:bvp} for the standard map \cref{eq:standard-map} as a function of the parameter $k$ for a product squared exponential kernel \cref{eq:sine-kernel} with varying $\sigma$, $N=500$, and $\epsilon = 10^{-5}$. As $\sigma$ increases, the residual becomes smoother and larger for large $k$. 
    (b) The residual for $\sigma = 0.1$, along with the energy contributions from $\Einv$ \cref{eq:Einv}, $\Ebd$ \cref{eq:Ebd}, and $\EK$ \cref{eq:EK}. For small $k$, the residual is dominated by the smoothness energy. As the amount of chaos increases, the residual becomes dominated by the invariance violation.
    }
    \label{fig:eval-plot}
\end{figure}

Now, we consider the dependence of the residual $R$ on the parameter $k$ (\cref{fig:eval-plot} (a)). 
We run the simulation scans with the same parameters \cref{fig:chaos-functions}, except we allow $k$ to vary between $0$ and $2$ and $\sigma$ to take $5$ values between $0.05$ and $0.5$. 
The residual $R$ smoothly and monotonically increases with $k$ for the three widest kernels $\sigma \in \{0.1, 0.2, 0.5\}$. 
For $\sigma = 0.1$ and $0.2$, the plot has two apparent phases, with a flatter section near $k = 0.0$ and a steeper section after.
For $\sigma = 0.5$, the initial flat phase is not visible.
The narrower kernels $\sigma \in {0.05, 0.075}$ are less smooth and more flat with changes in $k$.

We further investigate the $\sigma = 0.1$ plot in \cref{fig:eval-plot} (b).
In this, we plot the residual $R$, as well as the component energies $\Einv$, $\Ebd$, and $\EK$.
Initially, we see that the residual is dominated by the regularizing smoothing term $\EK$.
Then, the invariance energy $\Einv$ begins to dominate, causing a dramatic increase in $R$.
Thus, the two phases of the dependence of $R$ on $k$ can be explained by the transition from a smoothness penalty to an invariance penalty.

The above findings have implications if one wants to use \cref{alg:IBVP} for an optimization, as the parameters $\epsilon$ and $\sigma$ can be tuned to the application. 
The width $\sigma$ gives an indication of how smooth we require the invariant functions to be.
If $\sigma$ is large, then we have stronger requirements on $\sympmap$. 
The regularization $\epsilon$ sets the level of the flat section.
Between the two parameters, we can design the residual to permit a certain range of invariant functions, while more sharply penalizing functions outside that range.

\section{Conclusions}
\label{sec:conclusions}
In this paper, we have introduced two algorithms for finding invariant functions. 
The first, \cref{alg:IBVP}, solves a least-squares boundary value problem on a torus. 
This algorithm attempts to find a label function that foliates an area between two boundaries. 
The algorithm also returns a residual, which serves as a measure of our belief in the invariance of the function.
The second algorithm, \cref{alg:IEVP}, uses an eigenvalue problem to find a smooth label function with some zero boundary conditions.
In this case, the eigenvalue gives the measure of belief in the invariance. 
We have investigated the behavior of these algorithms with respect to their inputs $N$, $K$, $\sigma$, and $\sympmap$ in the examples \cref{sec:stellarator-example,sec:standard-map}.

There are many ways in which this work could be improved upon in the future. 
First, both algorithms could be faster with improvements in the numerical linear algebra.
Throughout this paper, except for the Arnoldi eigenvalue solver, we have only used simple dense linear algebra routines to solve \eqref{eq:bvp} and \eqref{eq:rayleighquotient}.
However, there are many ways to speed up kernel matrix evaluations and factorization, including multiplicative Schwarz \cite{beatson2001,leborne2019} and structured matrix approximations \cite{beatson1992,faul2005}. 
The evaluation of the symplectic map is also easily parallelized, but how to parallelize the rest of the algorithms is not obvious. 

There are also things that could be done to control error. 
For instance, one could consider using adaptive sampling to better resolve islands for fewer evaluations of the map.
One could also consider using leave-one-out validation for hyperparameter tuning or to extend the methods to use Gaussian processes. 
This could potentially give error bounds on the invariance measures, allowing for a more careful analysis of the invariance.

From a physical and dynamical systems point of view, it would also be interesting to extend this method to higher dimensional examples. 
The algorithms have been written in a way that they could be applied to higher dimensional examples, but we have not considered doing so in this paper. 
One natural first system would be the 4D standard map in \cite{firmbach2023}, where the boundary value problem could be used to test whether 3D surfaces foliate a 4D volume. 
This would provide a natural extension to gyro-averaged dynamics in stellarators, rather than the simple magnetic field line techniques. 
Additionally, in the higher dimensional problem, one could extend the method to learn Poisson-commuting label functions, allowing for the implicit definition of higher dimensional invariant tori.




\section*{Acknowledgments}
We thank Misha Padidar for his comments about boundary conditions and Benjamin Faber and Aaron Bader for the example stellarator configuration. We also thank Jim Meiss, Jackson Kulik, and William Clark.

\bibliographystyle{siam}
\bibliography{bibliography}

\end{document}